\documentclass[a4paper, 10pt]{amsart}
\usepackage{a4wide}
\usepackage{amsthm,amsmath,amssymb,amsfonts}
\usepackage[latin1]{inputenc}
\usepackage{enumerate}
\usepackage{varioref}
\usepackage[tips,arrow,matrix,curve]{xy}

\newcommand{\C}{\mathbb{C}}
\newcommand{\Nat}{{\mathbb{N}}}
\newcommand{\cover}{\gtrdot}
\def\setsuchas#1#2{\left\{\,{#1}\,\vrule\,{#2}\,\right\}}
\newcommand{\set}[1]{{\{#1\}}}
\newcommand{\commify}{\mathfrak{\sigma}}
\newcommand{\noncommify}{\mathfrak{\sigma^+}}
\newcommand{\e}[1]{\mathbf{e}_{#1}}
\newcommand{\f}[1]{\mathbf{f}_{#1}}
\newcommand{\dominant}{\mathcal{D}}
\newcommand{\NC}{\mathfrak{N}}
\newcommand{\NCn}[1]{{\NC}_{#1}}
\newcommand{\COMM}{\mathfrak{C}}
\newcommand{\COMMn}[1]{{\COMM}_{#1}}

\newcommand{\tdeg}[1]{\left \lvert {#1} \right \rvert}

\theoremstyle{definition}
\newtheorem{definition}{{Definition}}[section]

\theoremstyle{remark}

\theoremstyle{plain}
\newtheorem{lemma}[definition]{{Lemma}}

\newtheorem{proposition}[definition]{Proposition}
\newtheorem{theorem}[definition]{{Theorem}}

\newcommand{\Ya}{
            \put(0,0){\line(1,0){5}}
            \multiput(0,0)(1,0){6}{\line(0,1){1}}
            \put(0,1){\line(1,0){5}}
            \multiput(0,1)(1,0){6}{\line(0,1){1}}
            \put(0,2){\line(1,0){5}}
            \multiput(0,2)(1,0){6}{\line(0,1){1}}
            \put(0,3){\line(1,0){5}}
            \put(0,-1){\line(1,0){2}}
            \multiput(0,-1)(1,0){3}{\line(0,1){1}}
            \put(0,-2){\line(1,0){1}}
            \multiput(0,-2)(1,0){2}{\line(0,1){1}}
            \put(0,-3){\line(1,0){1}}
            \multiput(0,-3)(1,0){2}{\line(0,1){1}}
}

\newcommand{\Yb}{
            \put(0,0){\line(1,0){5}}
            \multiput(0,0)(1,0){6}{\line(0,1){1}}
            \put(0,1){\line(1,0){5}}
            \multiput(0,1)(1,0){6}{\line(0,1){1}}
            \put(0,2){\line(1,0){5}}
            \multiput(0,2)(1,0){6}{\line(0,1){1}}
            \put(0,3){\line(1,0){5}}
            \put(0,-1){\line(1,0){3}}
            \multiput(0,-1)(1,0){4}{\line(0,1){1}}
            \put(0,-2){\line(1,0){2}}
            \multiput(0,-2)(1,0){3}{\line(0,1){1}}
            \put(0,-3){\line(1,0){1}}
            \multiput(0,-3)(1,0){2}{\line(0,1){1}}
            \put(0,-4){\line(1,0){1}}
            \multiput(0,-4)(1,0){2}{\line(0,1){1}}
            \multiput(0,-1)(1,0){3}{\line(1,1){1}}
            \multiput(0,0)(1,0){3}{\line(1,-1){1}}
}

\newcommand{\Yd}{
            \put(2,-1){\line(1,1){1}}
            \put(2,0){\line(1,-1){1}}
           \put(0,0){\line(1,0){5}}
            \multiput(0,0)(1,0){6}{\line(0,1){1}}
            \put(0,1){\line(1,0){5}}
            \multiput(0,1)(1,0){6}{\line(0,1){1}}
            \put(0,2){\line(1,0){5}}
            \multiput(0,2)(1,0){6}{\line(0,1){1}}
            \put(0,3){\line(1,0){5}}
            \put(0,-1){\line(1,0){3}}
            \multiput(0,-1)(1,0){4}{\line(0,1){1}}
            \put(0,-2){\line(1,0){1}}
            \multiput(0,-2)(1,0){2}{\line(0,1){1}}
            \put(0,-3){\line(1,0){1}}
            \multiput(0,-3)(1,0){2}{\line(0,1){1}}
}

\newcommand{\YYa}{
\begin{picture}(5,6)(0,-4)    \thicklines \Ya \end{picture}
}

\newcommand{\YYb}{
\begin{picture}(5,6)(0,-4)    \thicklines \Yb \end{picture}
}

\newcommand{\YYd}{
\begin{picture}(5,6)(0,-4)    \thicklines \Yd \end{picture}
}

\begin{document}

\title{A poset classifying non-commutative term orders}

\author{Jan Snellman}
\address{Department of Mathematics\\
Link\"oping University\\
SE-58183 Link\"oping\\
Sweden \\
\texttt{jasne@mai.liu.se}
}


%\revision{2}
%\received{4 Feb 2001}
%\revised{20 Apr 2001}
%\accepted{16 Apr 2001}
\maketitle

\begin{abstract}
  We study a poset \(\NC\) on the free monoid \(X^*\) on a countable
  alphabet \(X\).  This poset is determined by the fact that its total
  extensions are precisely the standard term orders on \(X^*\). We
  also investigate the poset classifying degree-compatible standard
  term orders, and the poset classifying sorted term orders. For the
  latter poset, we give a Galois coconnection with the Young lattice.
\end{abstract}

%\tableofcontents
%\sloppy

\begin{section}{Introduction}
  So-called \emph{strongly stable ideals} are much studied in
  commutative algebra because of their intimate connection with
  \emph{generic initial ideals} 
  \cite{GreStill:GIN,Green:gin,Ebud:View}, because 
  their r\^ole in elucidating Macaulays theorem on possible
  Hilbert functions \cite{Bigatti95,Bigatti:Betti,BigattiRobbiano97}
  and because their minimal free resolutions have a simple structure
  \cite{Eliahou:MinRes}. In brief, a monomial ideal \(I\) in a
  polynomial ring \(\C[x_1,\dots,x_n]\) is strongly stable\footnote{In
    the literature, it is more common to insist on the reverse order
    of the variables, thus \(x_j\) is replaced with \(x_{j-1}\). For
    our purposes (particularly since we will be dealing with
    infinitely many variables) our definition is more convenient.} if,
  whenever \(m=x_1^{a_1} \cdots x_n^{a_n} \in I\), then
  \(\frac{x_{j+1}}{x_j} m \in I\) for all \(j\) such that \(j < n\),
  \(a_j >0\). So, we should be able to replace any occurring \(x_j\)
  with \(x_{j+1}\), as long as \(x_{j+1} \in \set{x_1,\dots,x_n}\). 
  
  Now let \(V\) denote the vector space of linear forms in
  \(\C[x_1,\dots,x_n]\), and let \(G_n\) denote the general linear
  group on \(V\). Then \(G_n\) acts on \(V\), and also on
  \(\C[x_1,\dots,x_n]\simeq S(V)\), the symmetric algebra on \(V\).
  This action is by linear substitution of variables. The ideals fixed
  by the subgroup of \(G_n\) consisting of the diagonal matrices are
  precisely the monomial ideals, and the ideals fixed by the subgroup
  consisting of the upper triangular matrices are precisely the
  strongly stable monomial ideals.

  Since a monomial ideal in \(\C[x_1,\dots,x_n]\) correspond to a
  monoid ideal in \([x_1,\dots,x_n]\), the free abelian monoid on
  \(\set{x_1,\dots,x_n}\), and since monoid ideals corresponds to
  filters in \(([x_1,\dots,x_n],D)\), where \(D\) denotes the
  divisibility partial order, it is natural to ask the question: is
  there a partial order \(\COMMn{n}\) on \([x_1,\dots,x_n]\), such
  that its filters 
  are precisely the strongly stable monomial ideals? And if so, what
  are its properties?

  It is clear that there is such a poset: simply define it as the
  smallest poset containing all relations \(m \le tm\) and
  \(m \le \frac{x_{i+1}}{x_i} m\). What is more interesting are the
  following two results, proved in \cite{Snellman:newglue}:

  \begin{enumerate}[A)]
  \item Define a \emph{standard term order} on \([x_1,\dots,x_n]\) to
    be a total order \(>\) such that \(1 < x_1 < \cdots < x_n\), and
    such that \(m_1 < m_2 \implies tm_1 < tm_2\). Then \(\COMMn{n}\)
    is the intersection of all standard term orders.
    \item Define a map from \([x_1,\dots,x_n]\) to the set of Ferrers
      diagrams with at most \(n\) columns, as follows:
      \(x_1^{a_1} \cdots x_n{a_n}\) goes to the diagram with \(a_i\)
      rows of length \(i\), for \(1 \le i \le n\). Then multiplication
      with \(x_j\) corresponds to inserting an extra row of length
      \(j\), and multiplication with \(\frac{x_{j+1}}{x_j}\)
      corresponds to removing one row of length \(j\) and inserting
      one row of length \(j+1\), i.e. to the insertion of an extra
      box (this is illustrated in Figure~\vref{fig:Ybox}). Thus, the
      map is an isotone bijection with isotone inverse, 
      showing that \(\COMMn{n}\) is isomorphic to a sub-lattice of the
      Young lattice. 
\setlength{\unitlength}{0.4cm}
     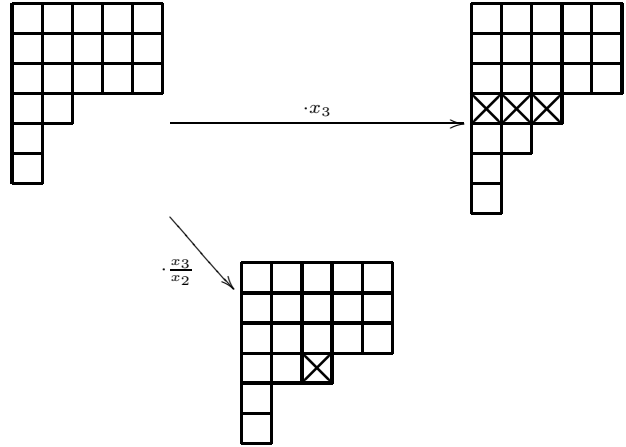
\begin{figure}[tbp]
        \begin{center}
      \begin{displaymath}
        \xymatrix{
          \YYa \ar@{->} [rr]^{\cdot x_3} 
          \ar@{->} [dr]_{\cdot \frac{x_3}{x_2}} && \YYb \\
          & \YYd
          } 
        \end{displaymath}
          \caption{\(x_1^2 x_2 x_5^3 \cdot x_3\) and \(x_1^2 x_2 x_5^3
            \cdot \frac{x_3}{x_2}\)} 
          \label{fig:Ybox}
        \end{center}
      \end{figure}
  \end{enumerate}
  
  Since the poset \(\COMMn{n}\) is canonically embedded in
  \(\COMMn{n+1}\), we can let \(n\) tend to infinity and study what
  happens when we have countably many indeterminates. It is a pleasing
  but unsurprising
  fact that the poset \(\COMM=\COMMn{\infty}\) that we obtain in this
  way is isomorphic to the Young lattice.

  In this article, we study the \emph{non-commutative} analogue
  \(\NCn{n}\) of the
  poset \(\COMMn{n}\). We show that this poset bears the same relation
  to non-commutative term orders, and to non-commutative monomial
  ideals fixed by upper triangular matrices, as does its commutative
  counterpart. However,   \(\NCn{n}\) is not a lattice, as we can see
  from Figure~\vref{fig:HasseNCN2}, the Hasse diagram of \(\NCn{2}\).

Our motivation for studying \(\NCn{n}\) is on one hand its relation to
``non-commutative generic initial ideals'', and on the other hand its
relation to term orders. We believe that it can be used to determine
the possible order types: recall the result of Martin and Scott
\cite{Martin:OT} that the possible order types for term orders on the
free monoid on two letters are \(\omega\), \(\omega^2\) and
\(\omega^\omega\). 

We also believe that  \(\NCn{n}\) is of interest ``in itself''; we
plan, in a subsequent article, to give an account of its incidence
algebra.

It is possible to define partial orders which capture the properties
of generic initial ideals over fields of finite
characteristic. Already in the commutative case,  
these
partial orders are immensely complicated; they involve the so-called
``Gauss order'' on the integers in a non-trivial way (see the PhD 
thesis of Keith
Pardue\cite{Pardue:Nonstandard} for more details.) Their
non-commutative counterparts should be even more formidable; we will
avoid these added complications and deal exclusively with the
characteristic zero case.
\end{section}

\begin{section}{Notations}
  We will use the terminology of \cite{Aigner:Combinatorial} for
  partially ordered sets.
Let \(n\) be a positive integer, 
  let \(X=\set{x_1,x_2,x_3,\dots}\) be a denumerable set of
  indeterminates, and let \(X_n =  \set{x_1,\dots,x_n}\). 
  Let \(X^*\), \(X_n^*\) be the corresponding free
  (non-abelian) monoid, and \([X]\), \([X_n]\) be the corresponding
  free abelian monoids. We let \(D\) denote the divisibility partial
  order on \(X^*\) and on \([X]\) (and, by restriction, on \(X_n^*\)
  and on \([X_n]\)). 

  \begin{definition}
    We define
    \begin{equation}
      \label{eq:commify}
      \begin{split}
      \commify: X^* & \to  [X] \\
      m=x_{i_1} \cdots x_{i_d} & \mapsto x_1^{a_1} \cdots x_N^{a_N}
      \end{split}
    \end{equation}
    where \(N\) is the highest index of a variable occurring in \(m\),
    and \(a_j\) denotes the number of occurrences of \(x_j\) in \(m\).

    We also define
    \begin{equation}
      \label{eq:noncommify}
      \begin{split}
      \noncommify: ([X],\COMM) & \to (X^*,\NC) \\
      x_1^{a_1} \cdots x_{\ell}^{a_\ell} & \mapsto x_1^{a_1} \cdots
      x_{\ell}^{a_\ell} 
      \end{split}
    \end{equation}
  \end{definition}

  \begin{lemma}\label{lemma:nccm}
   \(\commify\) and \(\noncommify\) are order-preserving
  with respect to the divisibility partial order on \(X^*\) and \([X]\).
    For \(m \in X^*\), \(t \in [X]\), we have that
    \(\commify(\noncommify(t))=t\) and that
    \(\noncommify(\commify(m))\) is the ``sorted version'' of \(m\);
    if \(m= x_{i_1} \cdots x_{i_d}\) then \(\noncommify(\commify(m)) =
    x_{i_{\tau(1)}} \cdots x_{i_{\tau(d)}}\) for some permutation
    \(\tau\) of \(\set{1,\dots,d}\).
  \end{lemma}
  \begin{proof}
    Obvious.
  \end{proof}

  Let \(\Nat^\omega\) be the subset of \(\Nat^{(\Nat^+)}\) consisting
  of finitely supported sequences, so \(\Nat^\omega\), with the
  natural order relation, is order isomorphic to \(\left( [X], D
  \right)\). We denote this isomorphism by \(\exp\), and its inverse
  by \(\log\). We put \(\e{i} = (0,\dots,0,1,0\dots)\), where the
  single 1 is in the \(i\)'th position.

  We let \(\Sigma\) be the functional which assigns to
  each element in \(\Nat^\omega\) the  sum of its components, and
  we let \(S\) be the (left) shift operator. The \(i\)'th projection
  map  \(\Nat^\omega \to \Nat\) is denoted by \(\pi_i\).
  
%  In what follows, we shall repeatedly make statements about
%  \(X^*\), for which there will be a corresponding statement about
%  \(X_n^*\). 
%  To avoid tedious repetitions, we will frequently use the phrase
%  ``the corresponding result holds for \(n\) variables'' to avoid duplicating
  
  \begin{definition}\label{def:multpart}
    Let \(M\) be a monoid. A partial order on \(M\) is
    \emph{multiplicative} (is a monoid partial order) if 
    \begin{itemize}
    \item \(1 < t\) for all \(t \in M \setminus \set{1}\),
    \item If \(s < t\) then  \(asb < atb\) for all   \(a,b \in M\).
    \end{itemize}
  \end{definition}

  \begin{definition}\label{def:termorder}
    A multiplicative total order \(<\) on \(X^*\)  (or on \(X_n^*\))
    is called a
    \emph{term order}. 
    A term order is \emph{standard} if
    \begin{equation}
      \label{eq:standard}
      x_1 < x_2 < x_3 < \cdots
    \end{equation}
  \end{definition}
  A result by Higman \cite{Higman:Lemma} implies that a standard term
  order is a well order. 
\end{section}

\begin{section}{Definition and basic properties of \protect\(\protect\NC\protect\)}
  \begin{definition}\label{def:R}
    Regarding a term order on \(X^*\) as a subset of
    \(X^* \times X^*\), we
    define the partial order \(\NC\) to be the intersection of 
    all standard term orders on \(X^*\). We define \(\NCn{n}\) to be
    the intersection of all standard term orders on \(X_n^*\).
  \end{definition}

  \begin{lemma}
    \(\NCn{n}\) is the restriction of \(\NC\).
  \end{lemma}
  \begin{proof}
    Every standard term order on \(X^*\) restricts to a standard term
    order on \(X_n^*\), and every standard term order on \(X_n^*\) can
    be extended to a standard term order on \(X^*\): just take the
    lexicographic product with some standard term order on \(\left(
      X \setminus X_n \right)^*\).
  \end{proof}

  \begin{proposition}
    \(\NC\) and \(\NCn{n}\) are locally finite well partial orders
    whose principal order ideals are finite.
  \end{proposition}
  \begin{proof}
    Let \(m \in X^*\). Since \(\NC\) is the intersection of all
    standard term orders, we have that if we pick one such standard
    term order, \(>\), then the principal order ideal on \(m\) with respect to
    \(\NC\) is contained in the principal order ideal on \(m\) with respect to
    \(>\). Since \(>\) is a well total order, this latter set is
    finite. 

    Any poset whose principal order ideals are finite is a locally
    finite well partial order, so the result follows.
  \end{proof}

  \begin{definition}\label{def:raising}
    For \(i \in \Nat^+\), define the
      \(i\)'th \emph{raising operation}  as the partially defined map
      \begin{equation}
        \label{eq:raise1}
        \begin{split}
        R_j : X^* & \to X^* \\
        m=x_{i_1} \cdots x_{i_d} & \mapsto 
          x_{i_1}\cdots x_{i_{j-1}}   x_{i_j +1} x_{i_{j+1}} \cdots
          x_{i_d} 
        \end{split}
      \end{equation}
      This is defined for \(j \le d\). As an operation from \(X_n^*
      \to X_n^*\), \(R_j(x_{i_1} \cdots x_{i_d})\) is defined when \(j
      \le d\) and \(i_j < n\).
  \end{definition}

  \begin{theorem}\label{thm:Rbasic}
        \begin{enumerate}[(i)] 
    \item \label{enum:text} The standard term orders on \(X^*\) 
      correspond to total multiplicative extensions of \(\NC\). Similarly, the
      standard term orders on \(X_n^*\) 
      correspond to total multiplicative extensions of \(\NCn{n}\).
    \item \label{enum:multiplic} \(\NC\) and \(\NCn{n}\) are monoid
      partial orders. 
    \item \label{enum:raising}  For \(m, m' \in X^*\),
      we have that \(m \le m'\) with respect to \(\NC\) if and only if 
      \(m'\) can be obtained from \(m\) by a finite sequence of
      applications of the following rules: 
      \begin{enumerate}[(a)]
      \item \label{enum:xm} \(t \mapsto x_1t\),
      \item \label{enum:mx} \(t \mapsto tx_1\),
      \item \label{enum:raise} \(t \mapsto R_j(t)\).
      \end{enumerate}
      The corresponding result holds for \(\NCn{n}\).
    \end{enumerate}   
  \end{theorem}
  \begin{proof}
    We prove the results for \(\NC\), the ones for \(\NCn{n}\) are similar.

    \eqref{enum:text} If \(<\) is a multiplicative total extension of
    \(\NC\), then it is a multiplicative total order on \(X^*\), hence
    a term order. Since it extends \(\NC\), it follows that \(x_1 <
    x_2 < x_3 < \cdots\) so \(<\) is standard.
    Conversely, if \(<\) is a standard term order, then \(x_1 <
    x_2 < x_3 < \cdots\) and \(<\) is a multiplicative total order on
    \(X^*\). Furthermore, if \(m \le m'\) with respect to \(\NC\),
    then \(m \le m'\) with respect to \emph{all} standard term orders,
    so in particular \(m < m'\) or \(m=m'\). Hence \(<\) extends \(\NC\).
    
    \eqref{enum:multiplic} Suppose that \(m \le m'\) with respect to
    \(\NC\). Then 
    \(m \le m'\) with respect to every standard term order, hence
    \(smt \le sm't\)  
    with respect to every standard term order, hence \(smt \le sm't\)
    with respect to \(\NC\). 

    \eqref{enum:raising} If \(>\) is a standard term order on \(\NC\), 
    then for all \(t\),
    \begin{displaymath}
      t < x_1t, \quad t < tx_1, \quad t \le R_j(t).
    \end{displaymath}
    Hence if \(m'\) is obtained from \(m\) by a sequence of operations 
    \eqref{enum:xm},     \eqref{enum:mx},     \eqref{enum:raise}, then 
    \(m' \ge m\). Since this holds for any standard term order, 
    \(m' \ge m\) with respect to \(\NC\).

    Conversely, suppose that \(m' \ge m\) with respect to
    \(\NC\). Then \(m' \ge m\) with respect to all standard term
    orders, in particular with respect to the total degree orders. So
    \(\tdeg{m'} \ge \tdeg{m}\). Furthermore, it is clear that regarded 
    as a subset of \(X^* \times X^*\), \(\NC\) is the smallest
    partial order which is also a bi-\(\NC\)-module containing
    \(\setsuchas{(smt, 
      m)}{s,t,m \in X^*}\) and \(\setsuchas{(x_j,x_i)}{j > i}\), the
    multiplication being  \(s(a,b)t = (sat,sbt)\). We can, in fact,
    replace these generators by the following: \(\setsuchas{(t,1)}{ t
      \in X^*}\), and \(\setsuchas{(x_{i+1}, x_i)}{i \in
      \Nat^+}\). But 
    \begin{align*}
      x_{i+1} &= R_1(x_i) \\
       x_{i_1}\cdots     x_{i_d} &= 
      R_d^{i_d} \circ \cdots \circ R_2^{i_2} \circ R_1^{i_1}(x_1^d) 
      \\
      &=
      R_d^{i_d} \circ \cdots \circ R_2^{i_2} \circ R_1^{i_1}(1\cdot
      x_1 \cdots x_1), 
    \end{align*}
    so \(x_{i+1}\) can be obtained from \(x_i\) by one application of
    a raising operator, and \(t=x_{i_1}\cdots  x_{i_d}\) can be
    obtained from \(1\) by a \(d\) right multiplications by \(x_1\),
    followed by an appropriate sequence of raising operators.
  \end{proof}

  Note that \(\NC\) and \(\NCn{n}\) have \emph{non-multiplicative
    total extensions} as 
  well. As an example, if we start extending \(\NCn{2}\) by first
  removing the anti-chain \(\set{x_1^2, x_2}\) by declaring that \(x_2
  > x_1^2\), then in order to have a multiplicative total extension, we must
  insist that \(x_1x_2 > x_1^3\), \(x_2x_1 > x_1^3\), \(x_2^2 >
  x_1^2x_2\), et cetera, and not the other way around. A
  non-multiplicative total extension may order these anti-chains
  independently.

  \begin{lemma}\label{lemma:covered}
    Let \(m = x_{i_1} \cdots x_{i_d}\), \(N=
    \max(\set{i_1,\dots,i_d})\). Let \(a_i\) denote the number of
    occurrences of \(x_i\) in \(m\); in other words,
    \((a_1,a_2,a_3,\dots) = \log(\commify(m))\). Then 
    \begin{enumerate}[(i)]
    \item 
      \(m\) is covered by the following words in \(X^*\): 
      \begin{itemize}
      \item \(x_1 m\) and \(m x_1\),
      \item The \(a_1\) words obtained by replacing one occurrence of
        \(x_1\) by \(x_2\), the \(a_2\) words obtained by replacing one
        occurrence of \(x_2\) by \(x_3\), and so on, up to and including
        the \(a_N\) words obtained by replacing one occurrence of
        \(x_N\) by \(x_{N+1}\). 
      \end{itemize}
      If \(m \neq x_1^k\), then these words are distinct, so that
      \(m\) is covered by exactly 
      \begin{displaymath}
        2 + \Sigma(\log(\commify(m))) = 2 + \sum_{i=1}^N a_i = 2 +
        \sum_{j=1}^d i_j 
      \end{displaymath}
      different words. On the
      other hand, \(x_1^k\)  
      is covered by \(x_{k+1}\), and by the \(k\) words \(x_1^a x_2
      x_1^{k-a-1}\), \(0 \le a \le k-1\).
      
    \item In  \(X_n^*\),  for \(N \le n\), we have that \(m\) is covered by
        \begin{itemize}
        \item \(x_1 m\) and \(m x_1\),
        \item The \(a_1\) words obtained by replacing one occurrence of
          \(x_1\) by \(x_2\), the \(a_2\) words obtained by replacing one
          occurrence of \(x_2\) by \(x_3\), and so on, up to and including
          the \(a_{n-1}\) words obtained by replacing one occurrence of
          \(x_{n-1}\) by \(x_n\). 
        \end{itemize}
        If \(m \neq x_1^k\), then these words are distinct, so that
        \(m\) is covered by exactly 
        \begin{displaymath}
           2 + \sum_{i=1}^{n-1} a_i 
        \end{displaymath}
        different words.

      \item
        The following words are covered by \(m\) (both in \(X^*\) and in 
        \(X_n^*\)):
        \begin{itemize}
        \item \(x_{i_2} \cdots x_{i_d}\), if \(i_1=1\),
        \item \(x_{i_1} \cdots x_{i_{d-1}}\), if \(i_d=1\),
        \item The \(a_2\) words obtained by replacing on occurrence of
          \(x_2\) with \(x_1\), and so on, up to and including
          the \(a_{N}\) words obtained by replacing one occurrence of
          \(x_{N}\) by \(x_{N-1}\). 
        \end{itemize}
        If \(m \neq x_1^k\), then these words are distinct, so that
        \(m\) covers exactly 
        \begin{displaymath}
          b + \sum_{i=2}^{n} a_i = b + \Sigma(S(\log(\commify(m)))), \qquad b = 
          \begin{cases}
            0 & i_1 \neq 1, \, i_d \neq 1 \\
            1 & i_1 \neq 1, \, i_d = 1 \\
            1 & i_1 = 1, \, i_d \neq 1 \\
            2 & i_1 = i_d = 1 \\
          \end{cases}
        \end{displaymath}
        different words. \(x_1^k\) covers exactly 1 word, namely
        \(x_1^{k-1}\). 
      \end{enumerate}
    \end{lemma}
    \begin{proof}
      This is immediate from Theorem~\vref{thm:Rbasic}.
    \end{proof}

        \setlength{\unitlength}{1.2cm}
  \begin{figure}[tbp]
    \begin{center}
      \begin{picture}(6,6)(0,0) \thicklines
\put(3.2, 0.0){\(1\)}
\put(3,0){\line(0,1){1}}

\put(3.1, 0.8){\(x_1\)}
\put(3,1){\line(-1,1){1}}
\put(3,1){\line(1,1){1}}

\put(1.4,1.9){\(x_1^2\)}
\put(2,2){\line(-1,1){1}}
\put(2,2){\line(1,1){1}}
\put(2,2){\line(2,1){2}}

\put(4.2,2){\(x_2\)}
\put(4,2){\line(-1,1){1}}
\put(4,2){\line(0,1){1}}

\put(0.5,2.9){\(x_1^3\)}
\put(2.4,2.9){\(x_1x_2\)}
\put(4.2,2.9){\(x_2x_1\)}

\put(1,3){\line(-1,1){1}}
\put(1,3){\line(1,1){1}}
\put(1,3){\line(2,1){2}}
\put(1,3){\line(3,1){3}}

\put(3,3){\line(-1,1){1}}
\put(3,3){\line(0,1){1}}
\put(3,3){\line(2,1){2}}

\put(4,3){\line(-1,1){1}}
\put(4,3){\line(0,1){1}}
\put(4,3){\line(1,1){1}}

\put(-0.5,4.1){\(x_1^4\)}
\put(1.5,4.1){\(x_1^2x_2\)}
\put(2.5,4.1){\(x_1x_2x_1\)}
\put(3.5,4.1){\(x_2x_1^2\)}
\put(5.1,4.1){\(x_2^2\)}

\end{picture}
      \caption{The Hasse diagram of \(\NCn{2}\).}
      \label{fig:HasseNCN2}
    \end{center}
  \end{figure}

 \end{section}

 \begin{section}{Strongly stable ideals}
   Let \(V\) be the complex vector space spanned by \(X\), and let
   \(G\) be the group of linear automorphisms of \(V\). Denote by
   \(T(V)\) the tensor algebra on \(V\). Then \(X^*\) is a basis of
   \(T(V)\) (as a vector space), and \(T(V) \simeq \C[X^*]\), the free
   non-commutative polynomial ring on \(X\). Furthermore, the action
   of \(G\) on \(V \simeq T(V)_1\) extends to an action on all of
   \(T(V)\): we define 
   \begin{equation}
     \label{eq:action}
     g.x_{i_1} \cdots x_{i_d} = (g.x_{i_1}) \cdots (g.x_{i_d})
   \end{equation}
   and extend this \(\C\)-linearly.

   By analogy with the commutative situation, we make the following
   definition: 

   \begin{definition}\label{def:upper}
     The subgroup of \emph{upper triangular} transformations in \(G\)
     is defined by 
     \begin{equation}
       \label{eq:upper}
       U = \setsuchas{u \in G}{u(x_i) =
         \sum_{j=i}^\infty c_{ij} x_j \text{ for all } i} 
     \end{equation}
   \end{definition}
   We note that the sums in \eqref{eq:upper} are finite, and that we
   must have that \(c_{ii} \neq 0\); otherwise, \(u\) would not be
   invertible.

   \begin{definition}\label{def:sstable}
     A monomial ideal in \(\C[X^*]\) is \emph{strongly stable} if it
     is fixed under the action of \(U\).
   \end{definition}
   
   \begin{theorem}\label{thm:sstable}
     The strongly stable monomial ideals in \(\C[X^*]\) correspond
     bijectively to filters in \(\NC\).
   \end{theorem}
   \begin{proof}
     Let \(I\) be a monomial ideal fixed by \(U\).
     Take any monomial \(m \in I\),  \(m=x_{a_1} \cdots c_{a_d}\).
     Define \(u \in U\) by \(u(x_{a_1})) = x_{a_1} + x_{a_1+1}\), \(u(x_j) =
     x_j\) for \(j \neq a_1\).  Then
   \begin{displaymath}
     u(m) = (x_{a_1} + x_{a_2}) u(x_{a_2} \cdots c_{a_d}) = m + m' +
     m'' + \cdots,
   \end{displaymath}
   where  \(m'=R_1(m)\), and the rest of the terms are similarly
   obtained from \(m\) by replacing one or several occurrences of
   \(x_{a_1}\) by \(x_{a_1+1}\). All those terms must be in \(I\),
   since \(I\) is a monomial ideal. By choosing different \(u\)'s, we
   get that \(I\) contains all monomials obtainable from \(m\) by
   means of a single raising operation. Since it is a monomial ideal,
   it contains also \(x_1m\) and \(m x_1\). Hence, from
   Theorem~\vref{thm:Rbasic} we get that the set of monomials in \(I\)
   is a filter with respect to \(\NC\).
   \end{proof}
   
   We get the corresponding result for the case of \(n\) variables:  we
   let \(V_n\) be the vector space spanned by \(X_n\), then \(T(V_n)
   \simeq \C[X_n^*]\), \(G_n\) is the general linear group on \(V\)
   and can be identified with the set of invertible \(n\times n\)
   matrices, and \(U_n\) with the set of upper triangular
   matrices. The \(n\) variable version of Theorem~\vref{thm:sstable}
   holds true.
 \end{section}

\begin{section}{The multi-ranking on \protect\(\protect\NC\protect\)}
Recall that a locally finite poset \((P, \ge)\) is \emph{ranked} if
there exists a \emph{rank function} \(\Phi: P \to \Nat\) such that
if \(m\) covers \(m'\) in \(P\), then \(\Phi(m) = \Phi(m')+1\). In
complete analogy, we define:
  \begin{definition}\label{def:multiranked}
    A locally finite poset \((P, \ge)\) is said to be
    \emph{\(\omega\)-multi-ranked} 
    if there exists a map 
    \begin{equation}
      \label{eq:multirankdef}
      \Phi: P \to \Nat^\omega
    \end{equation}
    such that 
    \begin{equation}
      \label{eq:covpfi}
      m \cover m' \quad \implies \quad \Phi(m) \cover \Phi(m')
    \end{equation}
    The poset is \emph{\(n\)-multi-ranked} if there exists   
    a map 
      \begin{equation}
        \label{eq:multirankdefn}
        \Phi_n: P \, \to \Nat^n
      \end{equation}
      such that 
      \begin{equation}
        \label{eq:covpfin}
        m \cover m' \quad \implies \quad \Phi_n(m) \cover \Phi_n(m')
      \end{equation}
  \end{definition}

  \begin{lemma}
    Let \(P\) be a locally finite poset, and let \(1 \le a \le
    b\). Then
    \begin{displaymath}
      P \text{ is \(\omega\)-ranked} \implies 
      P \text{ is \(b\)-ranked} \implies 
      P \text{ is \(a\)-ranked}  \implies
      P \text{ is \(1\)-ranked} \iff 
      P \text{ is ranked}.
    \end{displaymath}
  \end{lemma}
  
  As an example, we see that the Young lattice is \(\omega\)-ranked,
  with the multi-rank-function given by the natural embedding into
  \(\Nat^\omega\): thus a partition is multi-ranked by the sequence of
  lengths of the rows in its Ferrers diagram. Collapsing this ranking,
  we get the ordinary rank 
  function, which ranks an element in the Young lattice by the number
  of boxes in its Ferrers diagram.

  \begin{theorem}\label{thm:Rmranked}
    \(\NC\) is \(\omega\)-multi-ranked, and
      \(\NCn{n}\) is 
      \(n\)-multi-ranked.
  \end{theorem}
  \begin{proof}
   Let \(m=x_{i_1} \cdots x_{i_d} \in X^*\), and
    put 
    \[\mathbf{a} = \log(\commify(m))= (a_1,a_2,a_3,\dots) \in \Nat^\omega
    ,\]
    where \(a_j\) denotes the number of occurrences of \(x_j\) in \(m\).
    We give \(m\) multi-rank in the following way: 
    \begin{equation}
      \label{eq:mrank}
      \pi_j(\Phi(m)) = \Sigma(S^{j-1}(\mathbf{a}))
    \end{equation}
    We can decompose \(\Phi=G \circ \log \circ \commify\), where \(G\)
    is  
    the linear map 
    \begin{equation}
      \label{eq:mrank2}
      \begin{split}
        G: \Nat^\omega & \to \Nat^\omega \\
        G(\e{i}) & = \f{i} = \sum_{j=1}^i \e{j} 
      \end{split}
    \end{equation}

    Now suppose that \(m'\) covers \(m\). By
    Lemma~\vref{lemma:covered}, there are two cases:
    \begin{itemize}
    \item \(m'=x_1 m\) or \(m'=m x_1\). We see that 
      \begin{equation}
        \label{eq:diff1}
        \begin{split}
        \Phi(m')-\Phi(m) &= (G(\mathbf{a} + \e{1})) - G(\mathbf{a}) 
        \\
        &=
        (G(\mathbf{a}) - G(\e{1})) -  G(\mathbf{a})  \\ 
        &= G(\e{1})  \\ &= \e{1}          
        \end{split}
      \end{equation}
    \item \(m'\) is obtained from \(m\) by replacing
      one occurrence of \(x_j\) with an \(x_{j+1}\).  Then
      \begin{equation}
        \label{eq:diffj}
        \begin{split}
        \Phi(m')-\Phi(m) &= G(\mathbf{a} - \e{j} + e{j+1})
        -G(\mathbf{a}) \\
        & = -G(\e{j}) + G(\e{j+1}) \\
        &=  \e{j+1}          
        \end{split}
      \end{equation}
    \end{itemize}
    This shows that \(\Phi\) is a multi-rank function. The function
    \(\Phi_n\) is defined by restriction.
  \end{proof}

  By collapsing the ranking, we get
  \begin{theorem}\label{thm:rankgen}
    \(X^*\) and \(X_n^*\) are ranked posets. The rank of the word
    \(x_{i_1} \cdots x_{i_d}\), with \(a_j\) occurrences of the letter
    \(x_j\), is \(\sum_{j=1}^\infty j a_j\). The rank multi-generating
    functions for \(X^*\) and \(X_n^*\) are, respectively 
    \begin{equation}
      \label{eq:mrankgen}
      \frac{1}{1- \sum_{i=1}^\infty t^i x_i}, \quad \text{ and } \quad
      \frac{1}{1- \sum_{i=1}^n t^i x_i},
    \end{equation}
    the rank-generating functions are
    \begin{equation}
      \label{eq:rankgen}
      \frac{1-t}{1- 2t}, \quad \text{ and } \quad
      \frac{1-t}{1- 2t + t^{n+1}}
    \end{equation}
  \end{theorem}
  \begin{proof}
    We prove the formul{\ae}s for \(X^*\). The ranking \(\Phi\) gives
    \(x_i\)  weight \((1,\dots,1,0,\dots)\), with \(i\) consecutive
    ones. For  the standard weights, i.e. \(x_i\) has weight
    \((0,\dots,1,0\dots)\), with the 1 at position \(i\), the
    generating function is \(\frac{1}{1 - \sum_{i=1}^\infty
      x_i}\). Thus, substituting \(\prod_{j=1}^i t_j\) for \(x_i\), we
      obtain the rank-generating function 
      \begin{displaymath}
        \frac{1}{1 - \sum_{i=1}^\infty \prod_{j=1}^i t_j},
      \end{displaymath}
      which specialises to 
      \begin{displaymath}
        \frac{1}{1 - \sum_{i=1}^\infty \prod_{j=1}^i t} =
        \frac{1}{1 - \sum_{i=1}^\infty  t^i} = 
        \frac{1}{1 - t \frac{1}{1-t}} = \frac{1-t}{1- 2t}.
      \end{displaymath}
  \end{proof}

  In Figure~\vref{fig:HasseNC4}, we give the Hasse diagram of \(\NC\)
  up to and including rank level 3.
        \setlength{\unitlength}{1cm}
  \begin{figure}[htbp]
    \begin{center}
      \begin{picture}(12,6) 
        \put(5.9,0.7){\(1\)}
        \put(6,1){\line(0,1){1}}
        
        \put(5.3,2){\(x_1\)}
        \put(6,2){\line(-2,1){2}}
        \put(6,2){\line(2,1){2}}

        \put(3.1,2.6){\(x_1^2\)}
        \put(4,3){\line(-3,1){3}}
        \put(4,3){\line(1,1){1}}
        \put(4,3){\line(3,1){3}}

        \put(8.5,2.7){\(x_2\)}
        \put(8,3){\line(2,1){2}}
        \put(8,3){\line(-1,1){1}}
        \put(8,3){\line(-3,1){3}}

        \put(0.3, 4){\(x_1^3\)}
        \put(1,4){\line(-1,1){1}}
        \put(1,4){\line(1,1){1}}
        \put(1,4){\line(3,1){3}}
        \put(1,4){\line(5,1){5}}

        \put(4,4){\(x_1x_2\)}
        \put(5,4){\line(-1,1){1}}
        \put(5,4){\line(1,1){1}}
        \put(5,4){\line(2,1){2}}
        \put(5,4){\line(3,1){3}}

        \put(6,4){\(x_2x_1\)}
        \put(7,4){\line(-5,1){5}}
        \put(7,4){\line(-3,1){3}}
        \put(7,4){\line(0,1){1}}
        \put(7,4){\line(3,1){3}}

        \put(10.5,4){\(x_3\)}
        \put(10,4){\line(-2,1){2}}
        \put(10,4){\line(0,1){1}}
        \put(10,4){\line(2,1){2}}

        \put(0,5.3){\(x_1^4\)}

        \put(2,5.3){\(x_2x_1^2\)}

        \put(4,5.3){\(x_1x_2x_1\)}

        \put(6,5.3){\(x_1^2x_2\)}
        
        \put(7,5.3){\(x_2^2\)}
        
        \put(8,5.3){\(x_1x_3\)}

        \put(10,5.3){\(x_3x_1\)}

        \put(12,5.3){\(x_4\)}
      \end{picture}
      \caption{The Hasse diagram of \(\NC\).}
      \label{fig:HasseNC4}
    \end{center}
  \end{figure}
      
\end{section}

\begin{section}{Relation to commutative term orders}
  \begin{definition}\label{def:comm}
    Let \(\COMM\) denote the smallest partial order on the free
    abelian monoid \([X]\) such that 
    \begin{enumerate}
    \item \label{it:onesmall} 
      \(1 \le m\) for all \(m \in [X]\),
    \item \label{it:mult}
      \(m \le m' \implies tm \le tm'\) for all \(m,m',t \in [X]\),
    \item \label{it:repl}
      \(x_1^{a_1} \cdots x_n^{a_n} \le \frac{x_{i+1}}{x_i}
      x_1^{a_1} \cdots x_n^{a_n}\) whenever \(a_i > 0\).
    \end{enumerate}
    For any \(n\), denote by \(\COMMn{n}\) the restriction of this
    partial order to \([X_n] =[\set{x_1,\dots,x_n}]\).
  \end{definition}

  The following theorem was proved in \cite{Snellman:newglue} (parts 3
  and 4 can be easily derived from \cite{Pardue:Nonstandard}).
  \begin{theorem}\label{thm:comm}
    \begin{enumerate}
    \item \(\COMM\) is the intersection of all standard term orders on
      \([X]\). 
    \item \(\COMMn{n}\) is the intersection of all term orders on
      \([X_n]\). 
    \item     The map \(G \circ \log: ([X],\COMM) \to \Nat^\omega\) is an
    order-preserving  monoid monomorphism, and  \(([X],\COMM)\) is
    isomorphic to the image, which we call \(\dominant\).
    This image consists    of all 
    non-decreasing finitely supported sequences, and is
    order-isomorphic to the Young lattice of unordered
    number-partitions. 
  \item
    The image of \([X_n]\) correspond to the number-partitions whose
    diagrams have at most \(n\) columns.
    \end{enumerate}
  \end{theorem}
  In other words, \(G \circ \log\) is a \(\omega\)-multi-ranking,
  which happens to be an isomorphism onto its image.

  The relation between this poset and \(\NC\) is as follows.
  \begin{theorem}\label{thm:almostgalois}
    With respect to the partial order \(\NC\) on \(X^*\), and the
    partial order \(\COMM\) on \([X]\),
      \(\commify\) and \(\noncommify\) are monotone.
     Furthermore,   
     \(\noncommify(\commify(p)) \not < p\), 
     for all \(p     \in X^*\). 

     The same results hold for the restrictions of \(\commify\) and
     \(\noncommify\) to \(X_n^*\) and \([X_n]\).
  \end{theorem}
  \begin{proof}
    If 
    \begin{math}
X^*  \ni m =     x_{i_1} \cdots x_{i_d}
    \end{math}
     then 
     \begin{math}
       \commify( m) =
    x_1^{a_1} \cdots x_N^{a_N},
     \end{math}
 with \(a_\ell\) denoting
      the number  of \(j\) such that \(x_j = \ell\).
     Clearly,
    \(\commify(m x_1) = x_1     \commify(m)\), which is \(\ge
    \commify(m)\) with respect to \(\COMM\). 
    Furthermore, 
    \begin{displaymath}
      \begin{split}
\commify(R_j(m)) &= \commify(x_{i_1} \cdots
    x_{i_{j-1}} x_{i_j +1} x_{i_{j+1}} \cdots x_{i_d}) \\
    &= x_1^{a_1} \cdots
    x_{j-1}^{a_{j-1}} x_j^{a_j-1} x_{j+1}^{a_{j+1}+1}
    x_{j+2}^{a_{j+2}} \cdots x_N^{a_N} \\
    &= \frac{x_{j+1}}{x_j} m.
      \end{split}
    \end{displaymath}
    By the definition of \(\COMM\), this is \(\ge m\).

    Conversely, let \(t=x_1^{a_1} \cdots x_N^{a_N} \in [X]\). Then
    \(\noncommify(t) = x_1^{a_1} \cdots x_N^{a_N} \in X^*\). It is
    clear that \(\noncommify(x_1 t) = x_1 \noncommify(t) \ge
    \noncommify(t)\).  Finally,
    \begin{displaymath}
      \begin{split}
        \noncommify(\frac{x_{j+1}}{x_j} t) &= \noncommify( x_1^{a_1} \cdots
    x_{j-1}^{a_{j-1}} x_j^{a_j-1} x_{j+1}^{a_{j+1}+1}
    x_{j+2}^{a_{j+2}} \cdots x_N^{a_N}) \\
    &= x_1^{a_1} \cdots
    x_{j-1}^{a_{j-1}} x_j^{a_j-1} x_{j+1}^{a_{j+1}+1}
    x_{j+2}^{a_{j+2}} \cdots x_N^{a_N} \\
    &=R_{a_1 + a_2 + \cdots + a_{j}}(\noncommify(t))
      \end{split}
    \end{displaymath}
    So \(\noncommify\) is isotone.

    We recall (Lemma~\vref{lemma:nccm}) that
    \(\noncommify(\commify(p))\) is the ``sorted 
    version'' of \(p \in X^*\). Hence, \(p\) and
    \(\noncommify(\commify(p))\) have the same multi-rank, and form an
    anti-chain. 
  \end{proof}

The ranking \(\Phi\) is the composition of \(\commify\) and the
ranking of \(\COMM\), so the following diagram commutes:
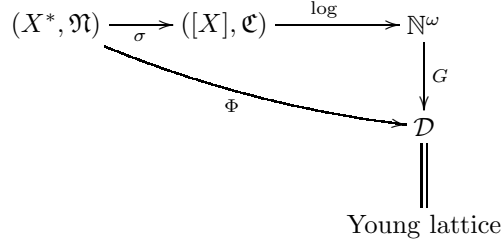
\begin{displaymath}
  \xymatrix{
    (X^*, \NC) \ar @/_/ [drr]_{\Phi}  \ar [r]_{\commify} &
    ([X], \COMM)  \ar [r]^{\log} 
    & \Nat^\omega \ar [d]^{G}\\ 
    && \dominant \ar @{=} [d] \\
    && \text{Young lattice}
    }
\end{displaymath}

We can thus regard
 \(\NC\) as a ``non-commutative version'' of the Young lattice.
A illustrative interpretation is the following: we identify the Young
lattice with \(\dominant\), and \(X^*,\NC\) with
all the set of all lattice walks in the infinite-dimensional lattice
\(\Nat^\omega\), using the steps \(\f{1} = \e{1} = (1,0,0,0,0\dots,)\),
\(\f{2} = \e{1}  + \e{2}= (1,1,0,0,0\dots,)\),
\(\f{3} = \e{1}  + \e{2} + \e{3}= (1,1,1,0,0,\dots,)\), et cetera. 
Then the multi-rank of such a walk is its endpoint, which lies in
\(\dominant\). When we restrict to 2 variables, the correspondence is
easy to draw, as is shown in Figure~\vref{fig:path}. It is easy to draw
the effect of the ``sorting'' \(\noncommify \circ \commify\), see
Figure~\vref{fig:pathc}. 

\begin{figure}[htbp]
  \begin{minipage}[t]{6cm}
        \setlength{\unitlength}{0.8cm}
\begin{picture}(6,6)
\multiput(0,0)(1,0){6}{\circle{0.2}}
\multiput(1,1)(1,0){5}{\circle{0.2}}
\multiput(2,2)(1,0){4}{\circle{0.2}}
\multiput(3,3)(1,0){3}{\circle{0.2}}
\multiput(4,4)(1,0){2}{\circle{0.2}}
\multiput(5,5)(1,0){1}{\circle{0.2}}

\put(0,0){\vector(1,1){1}}
\put(1,1){\vector(1,0){1}}
\put(2,1){\vector(1,0){1}}
\put(3,1){\vector(1,1){1}}
\put(4,2){\vector(1,1){1}}

\end{picture}  
    \caption{The word \(x_2 x_1^2 x_2^2\) as a path, bi-rank \((5,3)\)}
    \label{fig:path}
  \end{minipage}
 \hfill
  \begin{minipage}[t]{6cm}

        \setlength{\unitlength}{0.8cm}
\begin{picture}(6,6)
\multiput(0,0)(1,0){6}{\circle{0.2}}
\multiput(1,1)(1,0){5}{\circle{0.2}}
\multiput(2,2)(1,0){4}{\circle{0.2}}
\multiput(3,3)(1,0){3}{\circle{0.2}}
\multiput(4,4)(1,0){2}{\circle{0.2}}
\multiput(5,5)(1,0){1}{\circle{0.2}}

\put(0,0){\vector(1,0){1}}
\put(1,0){\vector(1,0){1}}
\put(2,0){\vector(1,1){1}}
\put(3,1){\vector(1,1){1}}
\put(4,2){\vector(1,1){1}}
\end{picture}    
    \caption{\(\noncommify(\commify(x_2 x_1^2 x_2^2))=x_1^2 x_2^3\)}
    \label{fig:pathc}
  \end{minipage}
\end{figure}

 \end{section}

\begin{section}{Variants}
  \begin{subsection}{``Sorted'' term orders}
Recall that if \(L,P\) are posets, 
then a
 \emph{Galois coconnection} is a pair of maps \(f: P \to L\)
  and \(g: L \to P\) such that
    \begin{enumerate}
    \item \(f\) and  \(g\) are
    order-preserving,
  \item \(gf (m) \ge m \) for all \(m \in P\),
  \item \(fg(w) \le w \) for all \(w \in L\).
    \end{enumerate}

  If we modify the partial order on \(X^*\) slightly,  \((\commify,
  \noncommify)\) becomes a  Galois  coconnection.

  \begin{definition}\label{def:sorted}
    Call a standard term order \(>\) \emph{sorted} if 
    \begin{equation}
      \label{eq:sorted}
        t x_i x_j s > t x_j x_i s \qquad \text{ for all } i < j, \,
        s,t \in X^*: 
    \end{equation}
 Let \(Q\) be the intersection of all sorted standard term orders on
    \(X^*\), and for each positive integer \(n\), \(Q\) be the
    intersection of all sorted standard term orders on 
    \(X_n^*\).
  \end{definition}
  
  \begin{lemma}\label{lemma:Qrest}
    For any commutative monomial \(m \in [X]\), the restriction of
    \(Q\) to \(\commify^{-1}(m)\) (equivalently, to the set of
    monomials having multi-rank \(G(\log(m))\)) is a chain.
  \end{lemma}
  \begin{proof}
    Trivial.
  \end{proof}

  \begin{theorem}\label{thm:galois}
    \(m \le m'\) with respect to \(Q\) if and only if \(m'\) can
    be obtained from \(m\) by a sequence of applications of the
    following types:

      \begin{enumerate}[(a)]
      \item \label{enum:nxm} \(t \mapsto x_1t\),
      \item \label{enum:nmx} \(t \mapsto tx_1\),
      \item \label{enum:nraise} \(t \mapsto R_j(t)\),
      \item \label{enum:nswap} \( t x_j x_i s \mapsto t x_i x_j s\)
        for \(i < j\), \(s,t \in X^*\).
      \end{enumerate}

    Furthermore,  \((\commify,
  \noncommify)\) is a  Galois  coconnection between \((X^*, Q)\) and
  \(([X], \COMM)\).

  If we restrict to \(n\) variables, we get a Galois coconnection
  between \((X_n^*, Q_n)\) and \(([X_n], \COMM)\).
  \end{theorem}
  \begin{proof}
    The first assertion is similar to Theorem~\vref{thm:Rbasic}; we
    omit the proof. 
    
    It remains to show that \((\commify,  \noncommify)\) is a Galois
    coconnection. 
    
     If \(m'\) is obtained from
    \(m\) by a sequence of operations of type
    \eqref{enum:nxm},\eqref{enum:nmx}, 
     \eqref{enum:nraise} then we know from
     Theorem~\vref{thm:almostgalois} that \(\commify(m') \ge
     \commify(m)\). Furthermore, \(\commify(t x_j x_i s) = \commify(t
     x_i x_j s)\). Hence, \(\commify\) is order preserving.
     
     Since \(Q \supset \NC\) ,  it follows from
     Theorem~\vref{thm:almostgalois} that  \(\noncommify\) is order preserving.

     For
     \(m=x_1^{a_1} \dots x_\ell^{a\ell}\) we
     have that \(\noncommify \commify (x_{i_1} \cdots x_{i_d}) =
     x_1^{a_1} \dots x_\ell^{a\ell}\), with \(a_j\) denoting the
     number of \(i\) such that \(a_i=j\). In other word, it is the
     ``sorted version'' 
     of \(m\), and it can be obtained from \(m\) by a sequence of
     operations of type \eqref{enum:nswap}: just perform a ``bubble
     sort''.  Hence \(\noncommify \commify (m) \ge m \) for all \(m \in
       X^*\).

       We know that t  \(\commify \noncommify(w) = w\) for all \(w \in [X]\),
       so in particular, 
       \(\commify \noncommify(w) \le w \) for all \(w \in
       [X]\).
  \end{proof}

  We note that with respect to the Galois coconnection above, the
  closed elements in \((X^*,Q)\) are the ``sorted'' 
  ones, and that all  elements in \([X], \COMM\) are coclosed.
  We show a part of the  Hasse diagrams for \(\COMMn{2})\) and
  \(Q_2\) in Figures~\vref{fig:hasseCom} and \vref{fig:hasseGal}.

\begin{figure}[htbp]
        \setlength{\unitlength}{0.45cm}
  \begin{minipage}[t]{7cm}

\begin{picture}(16,16)
\put(1,0){\circle*{0.2}}
\put(0,0){1}
\put(1,0){\line(0,1){4}}

\put(1,2){\circle*{0.2}}
\put(0,2){\(x_1\)}
\put(1,2){\line(0,1){4}}
\put(1,2){\line(2,1){4}}

\put(1,4){\circle*{0.2}}
\put(0,4){\(x_1^2\)}
\put(1,4){\line(0,1){4}}
\put(1,4){\line(1,1){4}}

\put(5,4){\circle*{0.2}}
\put(6,4){\(x_2\)}
\put(5,4){\line(0,1){4}}

\put(1,8){\circle*{0.2}}
\put(0,8){\(x_1^3\)}
\put(1,8){\line(0,1){6}}
\put(1,8){\line(2,3){4}}

\put(5,8){\circle*{0.2}}
\put(6,8){\(x_1x_2\)}
\put(5,8){\line(2,3){4}}
\put(5,8){\line(0,1){6}}

\put(1,14){\circle*{0.2}}
\put(0.5,15){\(x_1^4\)}

\put(5,14){\circle*{0.2}}
\put(5,15){\(x_1^2x_2\)}

\put(9,14){\circle*{0.2}}
\put(9,15){\(x_2^2\)}

\put(5,7.5){\oval(0.5,2)}
\put(5,13){\oval(0.5,3)}
\end{picture}    

    \caption{The Hasse diagram of \(\COMMn{2}\)}
    \label{fig:hasseCom}
  \end{minipage}
\hfill
  \begin{minipage}[t]{7cm}
\begin{picture}(10,16)
\put(1,0){\circle*{0.2}}
\put(0,0){1}
\put(1,0){\line(0,1){4}}

\put(1,2){\circle*{0.2}}
\put(0,2){\(x_1\)}
\put(1,2){\line(0,1){4}}
\put(1,2){\line(2,1){4}}

\put(1,4){\circle*{0.2}}
\put(0,4){\(x_1^2\)}
\put(1,4){\line(0,1){4}}
\put(1,4){\line(4,3){4}}

\put(5,4){\circle*{0.2}}
\put(6,4){\(x_2\)}
\put(5,4){\line(0,1){3}}

\put(1,8){\circle*{0.2}}
\put(0,8){\(x_1^3\)}
\put(1,8){\line(0,1){6}}
\put(1,8){\line(1,1){4}}

\put(5,7){\circle*{0.2}}
\put(6,7){\(x_2x_1\)}
\put(5,7){\line(0,1){1}}

\put(5,8){\circle*{0.2}}
\put(6,8){\(x_1x_2\)}
\put(5,8){\line(2,3){4}}
\put(5,8){\line(0,1){6}}

\put(1,14){\circle*{0.2}}
\put(0.5,15){\(x_1^4\)}

\put(5,14){\circle*{0.2}}
\put(5,15){\(x_1^2x_2\)}

\put(5,13){\circle*{0.2}}
\put(5.5,13){\(x_1x_2x_1\)}
\put(5,13){\line(0,1){1}}

\put(5,12){\circle*{0.2}}
\put(5.5,12){\(x_2x_1^2\)}
\put(5,12){\line(0,1){1}}

\put(9,14){\circle*{0.2}}
\put(9,15){\(x_2^2\)}

\put(5,7.5){\oval(0.5,2)}
\put(5,13){\oval(0.5,3)}

\end{picture}    
    \caption{The Hasse diagram of \(Q_2\) is obtained  from that of
      \(\COMMn{2}\) by replacing some elements by chains.}
    \label{fig:hasseGal}
  \end{minipage}

\end{figure}    

  \end{subsection}

  \begin{subsection}{Total degree term orders}
  \begin{theorem}\label{thm:degfix}
    Let \(n,d\) be non-negative integers. 
    \begin{enumerate}
    \item The restriction of \(\NC\) to the set of words in \(X^*\) of 
      total degree \(d\) is isomorphic to   \(\Nat^d\),
    \item The restriction of \(\NCn{n}\) to  the set of words in
      \(X_n^*\) of  
      total degree \(d\) is isomorphic to   the \(d\)-fold ordinal
      product of the chain with \(n\) elements,
    \item In particular, the restriction of \(\NCn{2}\) to the set of
      words in \(X_2^*\) of  
      total degree \(d\), is isomorphic to  the free boolean lattice
      on \(d\) elements.
    \end{enumerate}
  \end{theorem}
  \begin{proof} 
    The map 
\[X^* \ni x_{i_1} \cdots x_{i_d} \mapsto (i_1,\dots,i_d) \in \Nat^d\] 
is an    order-preserving bijection. If the word is in
\(X_n^*\) then \(1 \le i_j \le n\).
  \end{proof}

  \begin{theorem}\label{thm:degfirst}
    Let \(P\) denote the partial order on \(X^*\) which is the
    intersection of all \emph{degree-compatible} standard term orders,
    where a term order is degree-compatible if \(\tdeg{m} > \tdeg{m'}
    \implies m > m'\). Then \(P \subset \NC\), and \(P\) is the
    ordinal sum 
    \begin{displaymath}
      \Nat^0 + \Nat^1 + \Nat^2 + \cdots
    \end{displaymath}
    Similarly, if  \(P_n\) denotes the partial order on
    \(\set{x_1,\dots,x_n}^*\) which is the 
    intersection of all \emph{degree-compatible} standard term orders,
   then \(P_n \subset \NCn{n}\), and \(P_n\) is the
    ordinal sum 
    \begin{displaymath}
      C_n^0 + C_n^1 + C_n^2 + \cdots,
    \end{displaymath}
    where \(C_n\) is the chain with \(n\) elements. In particular,
    \(P_2\) is the ordinal sum 
    \begin{displaymath}
      B_0 + B_1 + B_2 + \cdots,
    \end{displaymath}
    where \(B_i\) is the free boolean lattice on \(i\) elements.
  \end{theorem}
  \begin{proof}
    Follows from Theorem~\vref{thm:degfix}.
  \end{proof}
    
  \end{subsection}

\end{section}

\bibliographystyle{plain}
%\bibliography{journals,articles,snellman}
\bibliography{nc}
\end{document}